\documentclass[11pt,oneside]{amsart}
\usepackage{amsmath,amssymb,graphics,array}
\newcolumntype{C}{>{$}c<{$}}

\newcommand{\bi}{\begin{itemize}}
\newcommand{\ei}{\end{itemize}}
\newcommand{\be}{\begin{enumerate}}
\newcommand{\ee}{\end{enumerate}}
\newcommand{\bc}{\begin{center}}
\newcommand{\ec}{\end{center}}
\newcommand{\bt}{\begin{tabular}}
\newcommand{\et}{\end{tabular}}

\newtheorem{thm}{Theorem}[section]
\newtheorem{lem}[thm]{Lemma}
\newtheorem{propn}[thm]{Proposition}

\theoremstyle{definition}

\newcommand{\blackboard}[1]{\ensuremath{\mathbb{#1}}}

\newcommand{\N}{\blackboard{N}}
\newcommand{\Z}{\blackboard{Z}}

\newcommand{\x}{\ensuremath \textbf{x}}

\newcommand{\ra}{\rightarrow}
\newcommand{\gm}{$\Gamma$}
\newcommand{\g}{\Gamma}
\newcommand{\pf}{poly-free}
\newcommand{\pfg}{poly-fg-free}
\newcommand{\rg}{right-angled Artin group}
\newcommand{\rgs}{right-angled Artin groups}

\newcommand{\sd}{\rtimes}
\newcommand{\ov}{\overline}

\newcommand{\nice}{doubly breakable cycle property}
\newcommand{\dd}{\delta}
\newcommand{\rn}{{\mathcal RN}}
\newcommand{\pfl}{{\mathsf{pfl}}}
\newcommand{\clq}{{\mathsf{clq}}}
\newcommand{\chr}{{\mathsf{chr}}}
\newcommand{\Aut}{{\mathsf{Aut}}}
\newcommand{\X}{{\mathbf{X}}}
\newcommand{\y}{\ensuremath \textbf{y}}
\newcommand{\z}{\ensuremath \textbf{z}}

\newcommand{\artin}{MR8:367a}
\newcommand{\bestvina}{MR2000h:20079}
\newcommand{\bestvbrady}{MR98i:20039}
\newcommand{\brieskorn}{MR54:10660}
\newcommand{\brown}{MR83k:20002}

\newcommand{\dicks}{MR82m:16024}
\newcommand{\droms}{MR88e:57003}
\newcommand{\feldmanhom}{MR0296168}
\newcommand{\green}{Egreen}
\newcommand{\hermeierartin}{MR2000m:20058}
\newcommand{\hermeiergrprod}{MR96a:20052}
\newcommand{\howiebb}{MR1713123}
\newcommand{\lyndonschupp}{MR58:28182}
\newcommand{\dmeierpf}{MR85h:20032}
\newcommand{\dmeierhomol}{MR82c:20091}
\newcommand{\meiervanwyk}{MR96h:20093}

\begin{document}
%%%%%%%%%%%%%%%%

\title[Poly-free constructions for right-angled Artin groups]
  {Poly-free constructions for right-angled Artin groups}

\author[S.~Hermiller]{Susan Hermiller$\!\,^1$}
\address{Dept. of Mathematics\\
        University of Nebraska\\
         Lincoln, NE 68588-0130}
\email{smh@math.unl.edu}

\author[Z.~Sunik]{Zoran $\check{\rm S}$uni$\acute {\rm k}$}
\address{Dept. of Mathematics\\
        Texas A\&M University\\
        College Station, TX 77843-3368}
\email{sunik@math.tamu.edu}

\begin{abstract}
We show that every right-angled Artin group $A\Gamma$
defined by a graph $\Gamma$ of finite chromatic number is
poly-free with poly-free
length bounded between the clique number and the chromatic number
of $\Gamma$. Further, a characterization of all right-angled Artin
groups of poly-free length 2 is given, namely the group $A\Gamma$ has
poly-free length 2 if and only if there exists an independent set
of vertices $D$ in $\Gamma$ such that every cycle in $\Gamma$
meets $D$ at least twice. Finally, it is shown that $A\Gamma$ is a
semidirect product of 2 free groups of finite rank if and only if
$\Gamma$ is a finite tree or a finite complete bipartite graph.
All of the proofs of the existence of \pf\ structures are
constructive.
\end{abstract}

\date{\today}
\maketitle

\footnotetext[1]{Supported under NSF grant no.\ DMS-0071037}

%%%%%%%%%%%%%%%%%%%%%%%%%%%%%%%%%%%%%%%%%%%%%%%%%%%%%%%%%%%%%%%%%%%%
\section{Introduction}\label{sec:intro}
%%%%%%%%%%%%%%%%%%%%%%%%%%%%%%%%%%%%%%%%%%%%%%%%%%%%%%%%%%%%%%%%%%%%

A group $G$ is {\it poly-free} if there exists a finite tower of
subgroups
$$
1 = G_0 \trianglelefteq G_1  \trianglelefteq \cdots \trianglelefteq G_N=G
$$
for which each quotient $G_{i+1}/G_i$ is a free group. The least
natural number $N$ for which such a tower exists is the {\it
poly-free length} of $G$, denoted $\pfl(G)$. A group $G$ is
poly-finitely-generated-free, or {\it poly-fg-free}, if there
exists a tower of this form with the additional property that each
of the quotients is a finitely generated free group. Since every
map onto a free group splits, a poly-free group can be realized as
an iterated semidirect product of free groups \cite{\dmeierpf}.

Examples of poly-fg-free groups include certain subgroups of Artin
groups. Let \gm\  be a finite simplicial graph; throughout the
text we will assume that such graphs do not have loops or multiple
edges.
If the edges of \gm\ are labeled by integers greater than one,
the associated {\it Artin group} $A\Gamma$ has generators
corresponding to the vertices, and relations
$$
{\underbrace{aba \cdots}_{n{\rm\;letters}}} = {\underbrace{bab
\cdots}_{n{\rm\;letters}}}
$$
where $\{a,b\}$ is an edge of the graph labeled $n$. If, in
addition, relations are added making each generator of order 2,
the resulting quotient is a Coxeter group.  {\it Braid groups} are
the Artin groups whose Coxeter quotients are the symmetric groups.
When the Coxeter quotient is finite, the Artin group is said to be
of {\it finite type}. {\it Pure} Artin groups are subgroups of
Artin groups which are the kernel of the homomorphism onto the
corresponding Coxeter group.

Pure braid groups are examples of poly-fg-free groups
\cite{\artin}, as are pure finite type Artin groups whose Coxeter
quotients are of type $B_n$, $D_n$, $I_2(p)$, and $F_4$
\cite{\brieskorn}. If the graph associated to an Artin group is a
tree, Hermiller and Meier \cite{\hermeierartin} have shown that
the Artin group is an extension of a finitely generated free group
by the integers, and hence is poly-fg-free. Recently, Bestvina
\cite{\bestvina}\  has asked if all Artin groups of finite type,
or indeed all Artin groups of any type, are virtually poly-free.

In this paper we investigate the poly-free properties of the class
of {\it right-angled Artin groups}, which are the Artin groups for
which the defining graph has every edge labeled 2.  That is, for a
graph \gm, the \rg\  $A\g$ is the group with
generators in one-to-one correspondence with the set $V(\g)$ of
vertices of \gm, and relations $[v,w]=vwv^{-1}w^{-1}$, for each
edge between vertices $v$ and $w$ of \gm.  These groups are also
known in the literature as graph groups, or free partially
commutative groups. (See \cite{\droms}, \cite{\green},
\cite{\hermeiergrprod}, \cite{\meiervanwyk}\  for information on
normal forms for right-angled Artin groups and further
references.)

Our main results are as follows.

\renewcommand{\thethm}{\Alph{thm}}
\begin{thm}\label{thm:rgispf}
Let \gm\  be a finite graph or, more generally, a
graph of finite chromatic number $\chr(\g)$ and finite
clique number $\clq(\g)$. The \rg\  $A\g$ is \pf. Moreover,
\[ \clq(\g) \le \pfl(A\g) \le \chr(\g), \]
and there exists a \pf\ tower for $A\g$ of length $\chr(\g)$.
\end{thm}

During the preparation of the text, W.~Dicks has pointed out to us
that J.~Howie~\cite{\howiebb} has established $|V(\g)|$ as an upper
bound for the poly-free length of a right-angled Artin group
defined by a finite graph $\Gamma$. Thus the above result is
an improvement in the case of finite graphs and a generalization
to a class of infinite graphs.

\begin{figure}[!hbt]
\begin{center}
\includegraphics{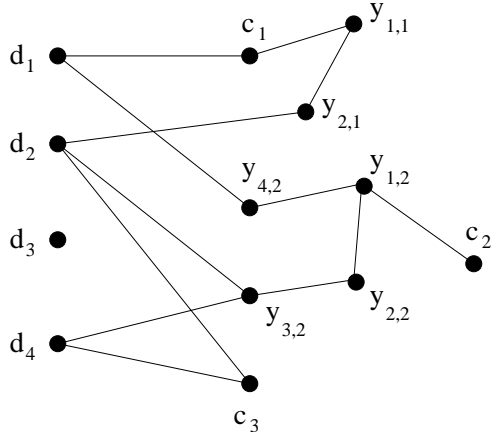}
\end{center}
\caption{An example of a graph with the doubly breakable cycle
property} \label{doubly}
\end{figure}

A graph \gm\  is said to have the {\it \nice} if \gm\ is not
totally disconnected and there exists a vertex subset $D \subseteq
V(\g)$ such that the full subgraph of \gm\  induced by $D$ is
totally disconnected, and such that every cycle in \gm\ contains
at least two vertices in $D$.
For a graph \gm\ with this property, the full subgraph
generated by the vertices in $V(\g)-D$ is a forest; if each of the
trees in this forest is collapsed to a point in $\g$, the
resulting graph is bipartite. Moreover, no vertex in $D$ is
connected by an edge to more than one vertex in each tree of the
forest $V(\Gamma)-D$. See Figure~\ref{doubly} for an example of a
graph with the doubly breakable cycle property. In this example
$D$ can be taken to be $D=\{d_1,d_2,d_3,d_4\}$.

\begin{thm}\label{thm:nicergispf2}
Let \gm\  be a graph. The \rg\  $A\g$ is \pf\  of
length 2 if and only if the graph \gm\  has the \nice.
\end{thm}

We note that
%both Theorem~\ref{thm:rgispf} and
Theorem~\ref{thm:nicergispf2} is valid in the context of
both finite and infinite graphs.

\begin{figure}[!hbt]
\begin{center}
\includegraphics{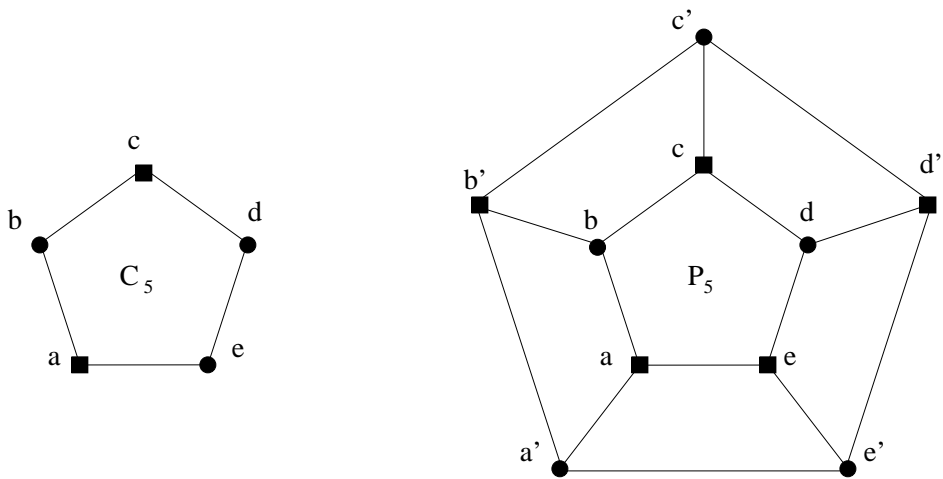}
\end{center}
\caption{A pentagon and a pentagonal prism}
\label{pentagon}
\end{figure}

Theorem~\ref{thm:nicergispf2} can be used to show that each of the
bounds in Theorem~\ref{thm:rgispf} are realized, as illustrated in
the following two examples. First, consider the graph $C_5$ given
by a 5-cycle (i.e. a pentagon; see Figure~\ref{pentagon}).  This
graph has chromatic number $\chr(C_5)=3$, so Theorem
\ref{thm:rgispf} shows that the group $AC_5$ is \pf\  with \pf\
length at most 3. However, the pentagon satisfies the \nice\ (for
example, one can take $D=\{a,c\}$). Thus Theorem
\ref{thm:nicergispf2} improves this bound to $\pfl(AC_5) \le 2$.
Indeed, since $\clq(C_5)=2$, this group contains $\Z^2$ as a
subgroup and is not free, so $\pfl(AC_5) = 2$. Hence the lower
bound on $\pfl(A\g)$ given by the clique number in
Theorem~\ref{thm:rgispf} is achieved in this example.
Next, suppose that $P_5$ is a pentagonal prism (see
Figure~\ref{pentagon}). In this case we have the same clique
number and the same chromatic number as for the pentagon, i.e.
$\clq(P_5)=2$ and $\chr(P_5)=3$, but $P_5$ does not satisfy the
\nice. Indeed, in a graph that satisfies the \nice\ exactly two
non-neighboring vertices must be selected from each 4-cycle to be
in the independent set of vertices $D$ breaking the cycles. Thus
if we choose $a$ in $D$, then we must also have $b'$, $c$, $d'$
and $e$ in $D$ (the vertices indicated by squares in the diagram
of $P_5$ in Figure~\ref{pentagon}). But $a$ and $e$ are neighbors,
so they cannot both be in $D$. This shows that $a$ cannot be in
$D$ and, by symmetry, no element can be in $D$. Since $D$ cannot
be empty $P_5$ does not satisfy the \nice.
Theorems~\ref{thm:rgispf} and \ref{thm:nicergispf2} show that
$\pfl(AP_5) \le 3$ and $\pfl(AP_5) \ne 2$. Thus $\pfl(AP_5) = 3$
and the chromatic number upper bound is achieved for this second
example.

If a graph \gm\ satisfies the \nice, the graph $\g$ can be colored
using three colors, one for the vertices in $D$, and two more for
the vertices in $V(\g)-D$ since their full subgraph is a forest.
Consequently, Theorem~\ref{thm:nicergispf2} implies that whenever
a \rg\ $A\g$ is \pf\ of length 2, then the defining graph \gm\
must have chromatic number at most three.

\begin{thm}\label{thm:pfgf2}
A right-angled Artin group $A\Gamma$ is poly-fg-free of length 2
if and only if $\Gamma$ is a finite tree or a finite complete
bipartite graph.
\end{thm}
\renewcommand{\thethm}{\thesection.\arabic{thm}}

For the right-angled Artin group $AC_5$ discussed above,
Theorem~\ref{thm:pfgf2} implies that the group $AC_5$ is not \pfg\
(of any length) even though it is poly-free of length 2. Indeed,
by the results of D.~Meier from~\cite{\dmeierhomol} (see
also~\cite{\feldmanhom}) the poly-fg-free length of a poly-fg-free
group is equal to its rational homological dimension. Since the
homological dimension of $AC_5$ is 2 (more on this later), $AC_5$
can only be poly-fg-free of length 2. However, $C_5$ is neither a tree
nor a complete bipartite graph and therefore $AC_5$ is not
poly-fg-free.

\subsection*{Organization}
Section~\ref{sec:background} is a brief review of \pf\ groups,
\rgs, and graph theoretic terminology. In Section~\ref{sec:pfimp}
we prove one direction of Theorem~\ref{thm:nicergispf2}, that
every \rg\ with \pf\ length 2 has the \nice, utilizing results of
\cite{\meiervanwyk} and \cite{\bestvbrady} on finiteness
properties of subgroups of \rgs. Section~\ref{sec:ineqs} contains
the proof of Theorem~\ref{thm:pfgf2}, utilizing a comparison of
Euler characteristics for \pf\ groups and \rgs, together with the
results of Section~\ref{sec:pfimp}. In Section~\ref{sec:rgpf} an
arbitrary \rg\ is exhibited as a split extension of a free
group by a \rg\ on a
subgraph, including an explicit description of the
action, proving Theorem~\ref{thm:rgispf}. Finally, in
Section~\ref{sec:niceimp2} we prove that for any graph
with the \nice, the corresponding \rg\ $A\g$ is a semidirect
product of two free groups, using a refinement of the techniques
of the previous section to construct the action. This result
completes the proof of Theorem~\ref{thm:nicergispf2}.

%%%%%%%%%%%%%%%%%%%%%%%%%%%%%%%%%%%%%%%%%%%%%%%%%%%%%%%%%%%%%%%%%%%%
\section{Background}\label{sec:background}
%%%%%%%%%%%%%%%%%%%%%%%%%%%%%%%%%%%%%%%%%%%%%%%%%%%%%%%%%%%%%%%%%%%%

\subsection{Groups}

Throughout the text, $g^a$ denotes the conjugate $a^{-1}ga$.

Let $G = \langle S \rangle$ be a group generated by $S$. A word
$w$ of length $k$ over $S \cup S^{-1}$ is a geodesic word if no
word over $S \cup S^{-1}$ of length strictly less than $k$
represents the same element in $G$ as $w$ does. A total order
defined on $S \cup S^{-1}$ induces a total order, called shortlex
order, on all words over $S \cup S^{-1}$ in which shorter words
always precede the longer ones and the words of the same length
are ordered lexicographically according to the order defined on $S
\cup S^{-1}$. A shortlex representative of an element $g \in G$ is the
smallest word in the shortlex order that represents $g$. Such a
representative is, by definition, geodesic.

\subsection{Graphs}

Throughout the paper, we assume that every graph is
a simplicial graph; that is, a simple undirected graph.
Therefore a graph $\Gamma$ is an ordered pair $\Gamma=(V,E)$ in
which the set $V=V(\Gamma)$ is a set of vertices and $E=E(\Gamma)$
is a set of edges, which is a set of two element subsets of $V$.
An edge $\{a,b\}$ has the vertices $a$ and $b$ as its endpoints.
Two vertices $x$ and $y$ in $V$ are neighbors (are adjacent) if
$\{x,y\}$ is an edge in $E$ (so no vertex is its own neighbor). A
cycle in $\Gamma$ is a path of length at least 3 in which no
vertex is repeated except for the initial and terminal one, which
coincide. The clique number $\clq(\Gamma)$ of a graph $\Gamma$ is
the largest size of a complete subgraph of $\Gamma$. Thus
$\clq(\Gamma)$ is the largest size of a subset $Q$ of $V$ for
which every 2-element subset $\{a,b\} \subseteq Q$ is an edge in
$E$. A proper coloring of a graph $\Gamma$ by $C$ is a labelling
$\ell:C \to V$ of the vertices in $V$ by symbols from a set of
colors $C$ in such a way that no two neighbors in $\Gamma$ are
colored in the same color. Thus if $\{a,b\} \in E$ then $\ell(a)
\neq \ell(b)$. The chromatic number of a graph $\Gamma$ is the
smallest size of a set $C$ for which there exists a proper
coloring of $\Gamma$ by $C$. A set of vertices $D$ is independent
if it can be colored by the same color in some proper coloring of
$\Gamma$. In other words, no two vertices in $D$ are adjacent.

\subsection{Right-angled Artin groups}

We freely use the following well known observation. If $\Gamma'$
is a subgraph of $\Gamma$ induced by a set of vertices $X
\subseteq V$, then the subgroup of $A\Gamma$ generated by the
elements of $X$ is $A\Gamma'$.

Throughout the text, given any homomorphism $\phi:A\g \to F_q$
from a \rg\ to a free group, the set $D :=\{v\in V(\g) ~|~
\phi(v)=1\}$ is called the set of {\it dead vertices}, the set
$L:=V(\g)-D$ is the set of {\it living vertices}, and the full
subgraph $\g_L$ generated by $L$ is the {\it living subgraph} of
\gm, with respect to $\phi$.

\begin{lem}\label{lem:raaggeod}\cite{\green}
Every geodesic representative of an element $t \in A\Gamma$ can be
obtained from any other representative of $t$ by finite number of
applications of the following operations:
\begin{enumerate}
 \item Eliminate a subword of the form $xx^{-1}$ or $x^{-1}x$
 with $x \in V(\g)$.
 \item If $x,y \in V(\g)$ are
 adjacent in $\g$, replace a single occurrence of
 $x^\pm y^\pm$ by $y^\pm x^\pm$.
\end{enumerate}
In particular, every geodesic
representative of an element $t \in A\Gamma_L$ can be obtained
from any other geodesic representative of $t$ by finite number of
applications of operation 2.
\end{lem}

\subsection{Poly-free groups}

Throughout the text, when $G$ is a semidirect product of two free
groups, we will write $G=F_k \rtimes F_q$ with associated
canonical homomorphism $\phi:G \to F_q$, so that the rank of the
kernel $\ker(\phi)$ is $k$ and the rank of the associate quotient
is $q$.

\begin{propn}\label{propn:subgroup}
If $G$ is \pf\  with length $N$ and $H \leq G$, then $H$ is \pf\
with length $\le N$.
\end{propn}

\begin{proof}
Given a \pf\   tower
$1 = G_0 \trianglelefteq G_1  \trianglelefteq \cdots \trianglelefteq G_N=G$
for $G$, then the tower
$1 = G_0 \cap H \trianglelefteq G_1\cap H  \trianglelefteq
\cdots \trianglelefteq G_N\cap H=H$ is a \pf\  tower for $H$.
\end{proof}

\begin{propn}\label{propn:quotient}
If $G$ has a normal free subgroup $H$ and the quotient $G/H$ is
\pf\  with $\pfl(G/H)=N$, then $G$ is \pf\  with $\pfl(G) \le
N+1$.
\end{propn}

\begin{proof}
Let $\phi:G \ra G/H$ be the canonical homomorphism.
Given a \pf\   tower
$1 = Q_0 \trianglelefteq Q_1  \trianglelefteq \cdots \trianglelefteq Q_N=G/H$
for $G/H$, then the tower
$1  \trianglelefteq \phi^{-1}(Q_0)=H \trianglelefteq \phi^{-1}(Q_1)
\trianglelefteq
\cdots \trianglelefteq \phi^{-1}(Q_N)=G$ is a \pf\  tower for $G$.
\end{proof}

%%%%%%%%%%%%%%%%%%%%%%%%%%%%%%%%%%%%%%%%%%%%%%%%%%%%%%%%%%%%%%%%%%%%
\section{Poly-freeness of length 2 implies the \nice}\label{sec:pfimp}
%%%%%%%%%%%%%%%%%%%%%%%%%%%%%%%%%%%%%%%%%%%%%%%%%%%%%%%%%%%%%%%%%%%%

Before proving the statement of the title of this section in
Proposition~\ref{propn:pf2impnicerg},
we begin with a few lemmas.

\begin{lem}\label{lem:concomp}
Let \gm\  be a graph with $A\g = F_k \rtimes F_q$ for free groups
$F_k$ and $F_q$ of finite or infinite rank.  Let $\phi:A\g \ra
F_q$ be the canonical homomorphism and let $\g_L$ be the
corresponding living subgraph. If $\g_\alpha$ is a connected
subgraph of $\g_L$, then the subgroup $\langle \{\phi(u) | u \in
\g_\alpha\}\rangle$ of $F_q$ is isomorphic to $\Z$.
\end{lem}

\begin{proof}
Fix a connected subgraph $\g_\alpha$ in $\Gamma_L$. The images
$\phi(u)$, $u\in \g_\alpha$ are nontrivial elements in the free
group $F_q$. Using the fact that $\g_\alpha$ is connected and
\cite[Prop. I.2.18]{\lyndonschupp} that says that the commuting
relation is an equivalence relation on the set of nontrivial
elements in a free group, we conclude that the group $\langle
\{\phi(u) | u \in \g_\alpha\}\rangle$ is abelian. However, the
only abelian subgroup of $F_q$ is $\Z$ and the conclusion follows.
\end{proof}

\begin{lem}\label{lem:adjdead}
Let \gm\ be a graph with $A\g = F_k \rtimes F_q$ for free groups
$F_k$ and $F_q$ of finite or infinite rank.  Let $\phi:A\g \ra
F_q$ be the canonical homomorphism and let $D :=\{v \in V(\g) ~|~
\phi(v)=1\}$ be the set of dead vertices.  If $d \in D$ and
$a_1,...,a_n \in L=V(\g)-D$ are all adjacent to $d$, then the
subgroup $\langle \phi(a_1),...,\phi(a_n) \rangle$ of $F_q$ is
free of rank $n$.
\end{lem}

\begin{proof}
Suppose that there exists a nontrivial word $u_1 \cdots u_m$ with
each $u_i \in \{a_1,...,a_n\}^{\pm 1}$ and $\phi(u_1) \cdots
\phi(u_m)=1$. Then $u_1 \cdots u_m \in \ker(\phi)$, and since each
$a_i$ is adjacent to $d$, the subgroup $\langle d,u_1 \cdots u_m
\rangle$ of $\ker(\phi)$ is isomorphic to $\Z^2$. This contradicts
the hypothesis that $\ker(\phi)=F_k$ is free.
\end{proof}

\begin{lem}\label{lem:coefind}
 Let \gm\  be a
finite graph, let $\phi:A\g \ra Z$ be any group homomorphism from
$A\g$ to an infinite cyclic group $Z=\langle z \rangle$, and let
$\rho:A\g \ra Z$ be the homomorphism defined by $\rho(v):=z$ for a
vertex $v\in V(\g)$ if $\phi(v) \ne 1$ and $\rho(v):=1$ if
$\phi(v)=1$. Then $\ker(\phi)$ is free if and only if $\ker(\rho)$
is free.
\end{lem}

\begin{proof}
Let $D:=\{w \in V(\g)~|~\phi(w)=1\}$ and for each $v \in V(\g)-D$,
let $n_v \in \Z$ be the unique integer such that $z^{n_v}
=\phi(v)$ in $Z$.

First suppose that the group $\ker(\rho)$ is not free. Let $N$ be
the least common multiple of the numbers $|n_v|$ for $v \in
V(\g)-D$.  Then for each $v \in V(\g)-D$, there exists an integer
$a_v$ such that $n_va_v=N$.  Define $a_w:=1$ for each $w \in D$.
Let $\mu:Z \ra Z$ be the homomorphism $z^k \mapsto z^{Nk}$ given
by taking the $N$-th power in $Z$. The composition $\mu \circ
\rho:A\g \ra Z$ has kernel $\ker(\mu \circ \rho)=\ker(\rho)$.
Define $\Theta:A\g \ra A\g$ by $\Theta(v):=v^{a_v}$ for all $v \in
V(\g)$; this defines a homomorphism of groups. Moreover, the
compositions $\phi \circ \Theta = \mu \circ \rho$. Hence
$\ker(\rho)=\ker(\mu \circ \rho) = \ker(\phi \circ \Theta)$.
$$\begin{matrix}
A\g &
\buildrel \Theta \over \longrightarrow &
A\g \\
\rho \downarrow&  &  \downarrow \phi \\
Z & \buildrel {\mu} \over \longrightarrow &
Z. \\
\end{matrix}$$

Suppose that $1 \ne g \in \ker(\Theta)$. Put a total ordering on
$V(\g)$, and let $g =_{A\g} v_1^{j_1} v_2^{j_2} \cdots v_k^{j_k}$
be the shortlex least representative of $g$ with each $v_i \in
V(\g)$, $j_i \in \Z$, and $v_i \ne v_{i+1}$. Then $\Theta(g)
=_{A\g} v_1^{j_1a_{v_1}} v_2^{j_2a_{v_2}} \cdots v_k^{j_ka_{v_k}}
=_{A\g} 1$.  Since $g \ne 1$, the word $v_1^{j_1} v_2^{j_2} \cdots
v_k^{j_k}$ is not empty, so there must be commutation relations
such that for some indices $i<m$ we have $v_i=v_m$, $j_ij_m <0$,
and $v_i$ commutes with $v_n$ for all $i <n<m$, so that
cancellation occurs between $v_i^{j_ia_{v_i}}$ and
$v_m^{j_ma_{v_m}}$. However, in that case cancellation must also
be possible in the normal form word $v_1^{j_1} v_2^{j_2} \cdots
v_k^{j_k}$, giving a contradiction. Therefore $\ker(\Theta) = 1$.

We now have that $\Theta$ is a monomorphism, and
$\Theta(\ker(\rho))= \Theta(\ker(\mu \circ \rho))=\Theta(\ker(\phi
\circ \Theta)) \le \ker(\phi)$.  Hence $\ker(\rho)$ is isomorphic
to a subgroup of $\ker(\phi)$.  Using the Nielsen-Schreier
Subgroup Theorem, that every subgroup of a free group is free,
 together with the hypothesis that $\ker(\rho)$ is
not free, we have that $\ker(\phi)$ cannot be a free group.

Next suppose that $\ker(\phi)$ is not free. Define the function
$\Psi:A\g \ra A\g$ by $\Psi(v):=v^{n_v}$ for $v \in V(\g)-D$ and
$\Psi(w):=w$ for $w \in D$.  Then $\phi=\rho \circ \Psi$.  An
argument similar to the proof above for $\Theta$ shows that $\Psi$
is a monomorphism of groups, so $\Psi$ restricts to an isomorphism
from $\ker(\phi)$ to a subgroup of $\ker(\rho)$. Since
$\ker(\phi)$ is not free, then $\ker(\rho)$ also cannot be free.
\end{proof}

Note that Lemma \ref{lem:coefind} is the
only lemma in this section for which the graph must be finite.
In the proof of Theorem B, this lemma is applied only
to a cycle $\Gamma'$ in $\Gamma$, which is finite, and
hence Theorem B is valid for both finite and infinite graphs $\Gamma$.

\begin{lem}\label{lem:z3}
The group $\Z^n$ has \pf\  length equal to $n$.
\end{lem}

\begin{proof}  Since $\Z^n= \underbrace{\Z \times \Z \times \dots \times \Z}_n$
is an iterated direct product of $n$ free groups, $\Z^n$ is \pf\
with length at most $n$.

To show that $\Z^n$ cannot have \pf\ length less than $n$, assume
that $\Z^n = (\dots (F_{n_1} \rtimes F_{n_2}) \rtimes \dots )
\rtimes F_{n_k}$, for some (finite or infinite) $n_i$,
$i=1,\dots,k$. Since $\Z^n$ is abelian, we must have
$n_1=n_2=\dots=n_k=1$. Thus $\Z^n$ can be generated by $k$
elements (one for each $F_{n_i}$). However, $Z^n$ cannot be
generated by fewer than $n$ elements, which shows that $n \leq k$.
\end{proof}

Note that Lemma \ref{lem:z3} shows that the \rg\  whose graph \gm\
is a triangle cannot be \pf\  of length 2.

\begin{propn}\label{propn:pf2impnicerg}
If \gm\  is a graph and the \rg\  $A\g$ is \pf\ of
length 2, then the graph \gm\  has the \nice.
\end{propn}

\begin{proof}
The group $A\g$ can be written as a semidirect product
$A\g=F_k \sd F_q$ where $F_k$ and $F_q$ are free groups
of ranks $k$ and $q$, respectively, and $k,q \in \N \cup \{\infty\}$.
Let $\phi:A\g \ra F_q$ be the canonical surjection,
with $\ker(\phi) =F_k$.
Let $D$ be the set of vertices $v$ in \gm\  with $\phi(v)=1$;
that is, $v \in ker(\phi)$.

If $d_1$ and $d_2$ are vertices in $D$, then the subgroup $\langle
d_1,d_2 \rangle < A\g$ generated by $d_1$ and $d_2$ must also be
contained in $\ker(\phi)$.  Since $\ker(\phi)$ is free, the
subgroup $\langle d_1,d_2 \rangle$ of $\ker(\phi)$ is free.  Hence
$d_1$ and $d_2$ cannot be joined by an edge in \gm. Thus the
subgraph induced by $D$ is totally disconnected, i.e., $D$ is an
independent set of vertices in $\Gamma$.

Suppose that $\g'$ is a cycle in the graph \gm. Then $\g'$
contains at least 3 vertices. The subgroup $A\g'$ of $A\g$ is
again a \rg, and Proposition~\ref{propn:subgroup} says that $A\g'$
is \pf\  with length at most 2. Lemma \ref{lem:z3} says that
$\Z^3$ is not \pf\  with length less than 3, so the cycle $\g'$
cannot have length 3. Thus $\g'$ contains at least 4 vertices.

If $V(\g') \cap D$ contains less than two vertices, then
the full subgraph of $\g'$ generated by the vertices in $V(\g')-D$
is connected.  In this case Lemma~\ref{lem:concomp}
shows that the restriction $\phi|_{A\g'}:A\g' \ra F_q$
has range $\langle f \rangle =\Z$ for some $f \in F_q$.
Since $\ker(\phi|_{A\g'})<\ker(\phi)$, then $\ker(\phi|_{A\g'})$
is also a free group.

If $V(\g') \cap D=\emptyset$, then $\phi(v) \ne f^0$ for all $v
\in V(\g')$.  Define the function $\rho:A\g' \ra \Z$ by
$\rho(v)=f$ for all $v \in V(\g')$.  Then Lemma \ref{lem:coefind}
says that $\ker(\rho)$ is free.  Since the graph $\g'$ is
connected, \cite[Theorem 6.3]{\meiervanwyk} shows that
$\ker(\rho)$ is finitely generated.  The cycle $\g'$ is a flag
complex, since $\g'$ is not a triangle, and this flag complex is
not simply connected.  The Main Theorem of \cite{\bestvbrady}
shows that $\ker(\rho)$ is not finitely presented.  Since a free
group cannot be finitely generated but not finitely presented, we
have a contradiction.

If there is exactly one vertex $d$ in $V(\g') \cap D$, then the
images $\phi(a)$ and $\phi(b)$ of the two neighbors $a$ and $b$ of
$d$ in the cycle $\Gamma'$ must generate a free group of rank 2 in
$F_q$ (by Lemma~\ref{lem:adjdead}). On the other hand, we already
established that the range of $\phi|_{A\g'}:A\g' \ra F_q$ is
cyclic, resulting again in a contradiction.
\end{proof}

%%%%%%%%%%%%%%%%%%%%%%%%%%%%%%%%%%%%%%%%%%%%%%%%%%%%%%%%%%%%%%%%%%%%
\section{Poly-fg-freeness of length 2}\label{sec:ineqs}
%%%%%%%%%%%%%%%%%%%%%%%%%%%%%%%%%%%%%%%%%%%%%%%%%%%%%%%%%%%%%%%%%%%%

In this section we prove Theorem~\ref{thm:pfgf2} using an
analysis of Euler characteristics of \rgs\ and \pf\ groups.

\begin{lem}\label{lem:reidschrei}
Let \gm\  be a finite graph with $A\g = F_k \rtimes F_q$
for free groups $F_k$ and $F_q$ of finite rank.  Let
$\phi:A\g \ra F_q$ be the canonical homomorphism, let
$D :=\{v \in V(\g) ~|~ \phi(v)=1\}$ be the set of
dead vertices, and for each $d \in D$, let $N_d$ denote
the set of vertices adjacent to $d$ in \gm.
Then the image of the subgroup generated by $N_d$
under the map $\phi$ has finite index in $F_q$, and
$k \geq  \sum_{d \in D} [F_q:\phi(\langle N_d\rangle)]~. $
\end{lem}

\begin{proof}
We begin by finding a presentation for the subgroup
$K:=\ker(\phi)=F_k$ of $A\g$ using the
Reidemeister-Schreier procedure, following the notation
in \cite[Proposition II.4.1]{\lyndonschupp}.
Let $F(V(\g))$ be the free group on the vertices of \gm,
let $\alpha:F(V(\g)) \ra A\g$ be the canonical
epimorphism, and let $\widetilde K:=\alpha^{-1}(K)$.
Then $F(V(\g))/\widetilde K \cong F_q$.
Choose a Schreier transversal $T$ for
$\widetilde K$ in $F(V(\g))$.

For any element $y \in F(V(\g))$, let $\ov{y}$ denote the
element of $T$ for which $Ky=K\ov{y}$.
For $t \in T$
and $v \in V(\g)$, define $\gamma(t,v):=tv\ov{tv}^{-1}$.
Given $d \in D$ and $t \in T$,
then  $\ov{td}=t$.  The element $d_t:=tdt^{-1}=\gamma(t,d)$ is a conjugate
of a nontrivial element in $A\g$, so $d_t$ itself is not trivial.
Given $a \in L$ and $t \in T$, define $a_t:=ta\overline{ta}^{-1} =
\gamma(t,a)$ as well.
The subset $S$ of nontrivial elements in the set
\[  \{ d_t \ | \ t \in T, \ d \in D \ \} \cup
   \{ a_t \ | \ t \in T, \ a \in L \ \} \]
generates $K$.

For $t \in T$ and $v \in V(\g)$, also define
$\gamma(t,v^{-1}):=tv^{-1}(\ov{tv^{-1}})^{-1}$. If $t \in T$, $d
\in D$, and $a \in L$, then $\gamma(t,d^{-1})=d_t^{-1}$ and
$\gamma(t,a^{-1})=(a_{\ov{ta^{-1}}})^{-1}$. Given any word $v=v_1
\cdots v_m$ with each $v_i \in V(\g)^{\pm 1}$, define
\[ \tau(v):=\gamma(1,v_1)\gamma(\ov{v_1},v_2) \cdots
 \gamma(\ov{v_1 \cdots v_{m-1}},v_m). \]
Note that for each $t \in T$ and each relator $r=[u,v] \in R$,
the element
$\tau(t)$ freely reduces to 1, and
\[ \tau(trt^{-1})=\gamma(t,u)\gamma(\ov{tu},v)\gamma(\ov{tuv},u^{-1})
 \gamma(\ov{tuvu^{-1}},v^{-1})~. \]
The latter words form the defining relators of the presentation
for $K=\ker(\phi)$. In particular, a defining set of relations is
given by
\begin{align*}
 R := & \{ d_t a_t = a_t d_{\overline{ta}} \ |
        \ t \in T, \ d \in D, \ a \in L,
        \ \{d,a\} \in E(\g) \ \}
        \ \cup \\
     & \{ a_t b_{\overline{ta}} = b_t a_{\overline{tb}} \ |
        \ t \in T, \ a,b \in L,
        \ \{a,b\} \in E(\g) \ \}
\end{align*}
and the group $K$ is presented by $\langle \ S \ | \ R \ \rangle$.

Abelianizing this presentation yields a presentation for
$K_{ab} = \Z^k$.  The subgroup $H$ of $K_{ab}$
generated by the elements of
\[ D_T := \{ \ d_t \ | \ d \in D, \ t \in T \ \} \]
is a free abelian direct factor of $K_{ab}$
presented by
\[
H =  \langle \ D_T \ | \ \{ d_t = d_{\overline{ta}} \ |
        \ t \in T, \ d \in D, \ a \in N_d \ \} \ \rangle_{ab}\ . \]
Since all relations in this presentation are equalities between
generators, the rank of $H$ is the number of equivalence classes
of generators. Note that two generators $d_t$ and $d_s$ are equal
in $H$ if and only if there exists a sequence of relations
$d_t=d_{\ov{ta_1}}=d_{\ov{ta_1a_2}}=\cdots=d_{\ov{ta_1 \cdots
a_m}}=d_s$ with $s=\ov{ta_1 \cdots a_{m+1}}$ and each $a_i \in
N_d^{\pm 1}$. This holds if and only if $Kta_1 \cdots a_{m+1}=Ks$,
which is satisfied if and only if $\phi(t)\phi(a_1) \cdots
\phi(a_{m+1})=\phi(s)$. Then $d_t$ and $d_s$ are equal in $H$ if
and only if $\phi(t)$ and $\phi(s)$ are in the same coset of
$\phi(\langle N_d \rangle)$ in $F_q$. Therefore the rank of $H$ is
equal to the sum of indices $\sum_{d \in D} [F_q:\phi(\langle N_d
\rangle)]$. Since the rank of $K_{ab}$ is $k$ we have
\[ k \geq rank(H) = \sum_{d \in D} [F_q:\phi(\langle N_d \rangle)]. \]
In particular, the
index $[F_q:\phi(\langle N_d \rangle)]$ is finite for all $d \in D$.
\end{proof}

\renewcommand{\thethm}{\ref{thm:pfgf2}}
\begin{thm}
A right-angled Artin group $A\Gamma$ is poly-fg-free of length 2
if and only if $\Gamma$ is either a finite tree or a finite
complete bipartite graph.
\end{thm}
\renewcommand{\thethm}{\thesection.\arabic{thm}}

\begin{proof}
If \gm\ is the complete bipartite graph $K_{k,q}$, then $A\g$ is
the direct product $F_k \times F_q$ of free groups of ranks $k$
and $q$.  On the other hand, if \gm\ is a tree on $n$ vertices,
the Artin group $A\g$ is a semidirect product $F_{n-1} \rtimes \Z$
\cite[Proposition 4.6]{\hermeierartin}.

Conversely, assume for the rest of this proof that $A\Gamma$ is a
poly-fg-free group of length 2. Since $A\Gamma$ is finitely
generated, the graph $\Gamma$ must be finite (since the
abelianization of $A\Gamma$ is $\Z^V$, any generating set of
$A\Gamma$ has at least $|V|$ elements).

There exists a split short exact sequence
\[ 1 \to F_k \to A\Gamma {\buildrel \phi \over \to} F_q \to 1, \]
such that the ranks $k$ and $q$ of the free groups are finite and
positive. By Proposition~\ref{propn:pf2impnicerg}, the graph \gm\
must satisfy the \nice.  Moreover, the proof of that proposition
shows that the set of dead vertices $D:=\{ v \in V(\Gamma) |
\phi(v)=1 \}$ associated to $\phi$ is an independent set of
vertices in $\Gamma$ such that every cycle in $\Gamma$ meets $D$
at least twice. Since $A\g$ is not poly-fg-free of length 1, the
living subgraph $\g_L$ is not empty.

Since semidirect products of free groups are torsion free, the Euler
characteristic of the semidirect product $A\g=F_k \rtimes F_q$ is
given by $\chi(A\Gamma) = \chi(F_k)\chi(F_q)$ (see
\cite[Proposition IX.7.3(d)]{\brown}). Therefore, given that the
Euler characteristic of a free group $F_r$ of rank $r$ is
$\chi(F_r)=1-r$,
\begin{equation} \chi(A\Gamma) = (k-1)(q-1)~. \label{eq:eulerkq}
\end{equation}

The Euler characteristic of a \rg\  can be computed using a
$K(A\g,1)$ space. The \nice\ implies that \gm\ contains edges but
does not contain any triangles. In this case the standard
2-complex associated to the standard presentation (from
Section~\ref{sec:intro}) of $A\g$ is a $K(A\g,1)$ \cite[Theorem
7.3]{\meiervanwyk}, and the associated chain complex is
\[ \cdots 0 \to \bigoplus_{|E(\g)|}\Z {\buildrel \partial_2 \over \to}
\bigoplus_{|V(\g)|}\Z {\buildrel \partial_1 \over \to} \Z \to 0~
\] (this complex can also be obtained as the co-invariants of the
free $\Z A\g$-module resolution in \cite{\dicks}). A
straightforward computation shows that all of the boundary maps
$\partial_i$ in this complex are trivial, so $H_0(A\g)=\Z$,
$H_1(A\g)=\Z^v$, and $H_2(A\g)=Z^e$ where $v:=|V(\g)|$ and
$e:=|E(\g)|$. Thus the Euler characteristic of $A\Gamma$ is also
\begin{equation} \chi(A\Gamma) = 1 - v + e~. \label{eq:eulerve}
\end{equation}

Denote the connected components of the nonempty graph $\Gamma_L$
by $C_1,\dots,C_c$.  Each of these components is a tree with
$n_j>0$ vertices, and hence $n_j-1$ edges, for $1 \le j \le c$.
Let $\dd:=|D|$ and denote the degrees of the vertices $d_1,\dots
d_\dd$ in $D$ by $g_1,\dots,g_\dd$. Rewriting
Equation~(\ref{eq:eulerve}) yields
\begin{multline}
 \chi(A\Gamma) =
1 - ( \dd + \sum_{j=1}^c n_j ) + ( \sum_{j=1}^c
 (n_j-1) + \sum_{i=1}^\dd g_i) =\\
 = 1 - \dd -c + \sum_{i=1}^\dd g_i = 1 -c + \sum_{i=1}^\dd (g_i-1).
\label{eq:eulercd}
\end{multline}

For each $1 \le i \le \dd$, let $N_i$ denote the set of vertices
adjacent to $d_i$. Lemma~\ref{lem:adjdead} says that the subgroup
$\phi(\langle N_i \rangle)$ of $F_q$ is free of rank equal to the
degree $g_i$ of $d_i$. Using the Schreier Formula,
$(rank(F_q)-1)[F_q:\phi(\langle N_i\rangle)]=
 (rank(\phi(\langle N_i\rangle)-1)$, so
\[ g_i-1 = (q-1)[F_q:\phi(\langle N_i\rangle)]. \]
According to Lemma~\ref{lem:reidschrei}, $\sum_{i=1}^\dd
[F_q:\phi(\langle N_{d_i}\rangle)] \le k$. Thus, taking into
account the non-negativity of $q-1$,
\[ \sum_{i=1}^\dd (g_i -1) = (q-1)\sum_{i=1}^\dd [F_q:\phi(\langle N_{d_i}\rangle)] \leq k(q-1). \]
Combining this with Equation~(\ref{eq:eulerkq}) and
Equation~(\ref{eq:eulercd}), then
\[(k-1)(q-1) ~~{\buildrel (\ref{eq:eulerkq}) \over =}~~
\chi(\g) ~~{\buildrel (\ref{eq:eulercd}) \over =}~~ 1 -c +
\sum_{i=1}^\dd (g_i-1) ~\leq~ 1 -c + k(q-1)~, \] which implies
that $c \leq q \ . $

Lemma \ref{lem:concomp} says that for each component $C_j$, there
exists an element $f_j$ in $F_q$ such that all vertices from the
component $C_j$ are mapped by $\phi$ to a power of $f_j$. Since
$\phi$ is onto and the dead vertices in $D$ are mapped to the
identity in $F_q$, $f_1,\dots,f_c$ generate $F_q$, which implies
that $ q \leq c.$ Using the inequality at the end of the previous
paragraph, then
\begin{equation} q = c. \label{eq:qec}
\end{equation}

If $q=c=1$, then the living subgraph $\g_L$ of \gm\  is a single
tree. Since the kernel $\ker(\phi)=F_k$ is finitely generated,
\cite[Theorem 6.1]{\meiervanwyk} says that
%$\Gamma$ is connected. Since $D$ is an independent set of vertices,
every dead vertex in
$D$ must be adjacent to a vertex in $\g_L$. The \nice\ says that
there cannot be a cycle in \gm\ that meets a dead vertex in $D$
only once, which implies that for each $d \in D$, $d$ cannot be
the endpoint of two different edges whose other endpoints lie in
$\g_L$. Therefore in this case the graph \gm\ is also a tree.

Finally suppose that $q=c \geq 2$. As in the previous paragraph,
the \nice\  says that for each $d_i \in D$, $d_i$ cannot be the
endpoint of two different edges whose other endpoints lie in the
same component $C_j$ of $\g_L$. Hence each degree $g_i \le c$, so
$\sum_{i=1}^\dd g_i \leq c\dd$. Using Equation~(\ref{eq:qec}),
Equation~(\ref{eq:eulerkq}), and Equation~(\ref{eq:eulercd}), then
\begin{eqnarray*}
(k-1)(c-1) & {\buildrel (\ref{eq:qec}) \over =} & (k-1)(q-1)
~~{\buildrel (\ref{eq:eulerkq}) \over =}~~ \chi(A\Gamma)
 ~~{\buildrel (\ref{eq:eulercd})  \over =}~~ 1 - \dd -c + \sum_{i=1}^\dd g_i\\
& \leq &  1 - \dd -c +c\dd = (c-1)(\dd-1)~.
\end{eqnarray*}
Since $c \ge 2$, therefore
\begin{equation} k \leq \dd~. \label{eq:kld}
\end{equation}
Note that  equality holds if and only if each vertex in $D$ has degree
$c$; i.e., there exists a single edge between each vertex in $D$ and
each component of $\Gamma_L$.

Since $A\Gamma = F_k \rtimes F_q$, the group $A\Gamma$ can be
generated by $k+q$ elements. However, the minimal number of
generators for the \rg\ $A\Gamma$
is the number $v=\dd + \sum_{j=1}^c n_j$
of vertices in $\Gamma$, which
in turn is at least as large as $\dd+c$. Thus
\begin{equation} \dd + c \leq k + q ~.
\label{eq:dclkq}
\end{equation}
Note that equality in this case is possible only
if $\dd+c=v=\dd + \sum_{j=1}^c n_j=k+q$, and
so $c=\sum_{j=1}^c n_j$.  Thus equality implies that each $n_j=1$;
i.e., each component $C_j$ of $\Gamma_L$ is a single vertex.

Using the fact that $q=c$, Inequality~(\ref{eq:kld}) and
Inequality~(\ref{eq:dclkq}) imply that $\dd=k$, so equality holds
both in (\ref{eq:kld}) and in (\ref{eq:dclkq}). Therefore $D$ and
$L$ are each independent sets of vertices in \gm, and there is an
edge between each vertex in $D$ and each vertex in $L$. Thus in
this case $\g$ is a complete bipartite graph.
\end{proof}

%%%%%%%%%%%%%%%%%%%%%%%%%%%%%%%%%%%%%%%%%%%%%%%%%%%%%%%%%%%%%%%%%%%%
\section{Every \rg\  is poly-free}\label{sec:rgpf}
%%%%%%%%%%%%%%%%%%%%%%%%%%%%%%%%%%%%%%%%%%%%%%%%%%%%%%%%%%%%%%%%%%%%

%In this section we prove the following.
Given a graph \gm\ with finite chromatic number greater
than one, let $D$ be the set of vertices in one of the colors and
let $L:=V(\g)-D$ be the vertices in the other colors.  Let $\g_L$
be the full subgraph of \gm\ induced by $L$. If we define a
homomorphism $\phi:A\g \ra A\g_L$ by $\phi(d)=1$ for $d \in D$ and
$\phi(a)=a$ for $a \in L$, then $D$ is the set of dead vertices
and $\g_L$ is the living subgraph associated to this homomorphism.
In the following
proof we construct a free group $F$ (isomorphic to $\ker(\phi)$)
and an action of $A\g_L$ on $F$, and exhibit directly that $A\g$
is isomorphic to the semidirect product $F \rtimes A\Gamma_L$ of a
free group with $A\g_L$.

\renewcommand{\thethm}{\ref{thm:rgispf}}
\begin{thm}
Let \gm\  be a finite graph or, more generally, a
graph of finite chromatic number $\chr(\g)$ and finite
clique number $\clq(\g)$. The \rg\  $A\g$ is \pf. Moreover,
\[ \clq(\g) \le \pfl(A\g) \le \chr(\g), \]
and there exists a \pf\ tower for $A\g$ of length $\chr(\g)$.
\end{thm}
\renewcommand{\thethm}{\thesection.\arabic{thm}}

\begin{proof}
To prove poly-freeness and the upper bound on the \pf\ length, we
induct on $\chr(\g)$.  If $\chr(\g)=1$, then $\g$ is totally disconnected,
so $A\g$ is free, and hence \pf\ of length 1.  Next suppose that
$\chr(\g) \geq 2$, and that
for every graph $\g'$ with $\chr(\g') < \chr(\g)$, the group
$A\g'$ is poly-free and has a \pf\ tower of length $\chr(\g')$.

Choose a coloring of $\Gamma$ in $\chr(\g)$ colors, one of which
is gray. Let $D$ be the set of vertices in $V=V(\Gamma)$
colored in gray, $L=V-D$ be the set of vertices colored in a
different color, $\Gamma_L$ be the subgraph of $\Gamma$
induced by $L$ and $A\Gamma_L$ be the corresponding \rg.
Then $\chr(\Gamma_L)=\chr(\Gamma)-1$ and the inductive
assumption implies that there exists a poly-free tower for
$A\Gamma_L$ of length $\chr(\Gamma)-1$.

In the discussion that follows, a geodesic representative of an
element $t \in A\Gamma_L$ means a geodesic word in the alphabet
$L^{\pm 1}$. For any vertex $v \in V(\g)$, denote by $N_v$ the set
of vertices adjacent to $v$; i.e., the neighbors of $v$. For each
$d \in D$, define a set of symbols
\[ T_d := \{\ d_t \ | \ t\in A\Gamma_L, \text{ no
  geodesic rep. of } t
  \text{ starts with a letter in } N_d^{\pm 1} \ \}. \]
Let $F(T_d)$ be the free group over $T_d$ and let $F$ be the free
group
\[ F := *_{d \in D} F(T_d). \]

For each generator $a \in L$,
define an endomorphism $\alpha_a:F \ra F$ by
\begin{equation}\label{AGL action}
\alpha_a(d_t) := \begin{cases} d_{ta}, & d_{ta} \in T_d \\
                         d_t,  & d_{ta} \not \in T_d
       \end{cases},
\end{equation}
for all $d \in D$ and $d_t \in T_d$. Since $F$ is a free group,
this definition of $\alpha_a$ on the generators of $F$ extends to
an endomorphism on $F$. In order to show that $a \mapsto \alpha_a$
extends to an action of $A\g_L$ on $F$, we first need to consider
when the conditions $d_t \in T_d$ and $d_{ta} \not \in T_d$ occur
simultaneously.

Assume that $d_t \in T_d$. If $d_{ta} \not \in T_d$, then there
exists a geodesic representative $w$ of $ta$ that begins with a
letter in $N_d^{\pm 1}$. Consider the word $wa^{-1}$ representing
$t$. Since $d_t \in T_d$ the word $wa^{-1}$ cannot be geodesic. By
Lemma~\ref{lem:raaggeod} we can write $w$ as $u_1au_2$ where $a$
commutes with all the letters in $u_2$. The word $u_1u_2$ is a
geodesic representative of $t$. As such, it cannot start in
$N_d^{\pm 1}$. Thus the geodesic word $w=u_1au_2$ representing
$ta$ and the geodesic word $u_1u_2$ representing $t$ start in a
different letter. This is possible only when $u_1$ is empty. Thus
$w=au_2$ is a geodesic representative of $ta$ that starts in
$N_d^{\pm 1}$ and $a$ commutes with all letters in $u_2$, i.e.,
$a$ commutes with $d$ and all the letters in $u_2$. However, $u_2$
is a geodesic representative of $t$. Since any other geodesic
representative of $t$ can be obtained from $u_2$ by commuting
letters we conclude that $a$ commutes with $d$ and all the letters
in any geodesic representative of $t$. A similar proof shows that,
conversely, if $a$ commutes with $d$ and all letters in any
geodesic representative of $t$ then $d_{ta}$ cannot be in $T_d$.
Therefore
\begin{description}
\item [$(*)$] If $d_t \in T_d$ then $d_{ta} \not \in T_d$   if and
only if $d$ and all of the symbols in any geodesic representative
of $t$ are adjacent to $a$ in $\Gamma$.
\end{description}

For each $a \in L$ and $d_t \in T_d$,
define another endomorphism $\alpha_{a^{-1}}:F \ra F$ by
replacing $a$ by $a^{-1}$ in Equation~\ref{AGL action}.  Then
\[ \alpha_a(\alpha_{a^{-1}}(d_t)) =
\alpha_a\left( \begin{cases} d_{ta^{-1}}, & d_{ta^{-1}} \in T_d \\
                         d_t,  & d_{ta^{-1}} \not \in T_d
       \end{cases} \right)  =
 \begin{cases} d_{t}, & d_{ta^{-1}} \in T_d, \ d_{ta} \in T_d \\
               d_{ta^{-1}},  & d_{ta^{-1}} \in T_d, \ d_{ta} \not \in T_d\\
               d_{ta},  & d_{ta^{-1}} \not \in T_d, \ d_{ta} \in T_d\\
               d_t,     & d_{ta^{-1}} \not \in T_d, \ d_{ta} \not \in T_d
 \end{cases}.
\]
As a consequence of $(*)$ from the previous paragraph,
for every $a \in L$, $d \in D$,
and $d_t \in T_d$, we have $d_{ta} \in T_d$ if and
only if $d_{ta^{-1}} \in T_d$, and so the middle two cases
in last expression of the equation above cannot occur.
Therefore $\alpha_a(\alpha_{a^{-1}}(d_t))=d_t$, and
similarly $\alpha_{a^{-1}}(\alpha_a(d_t))=d_t$.  Thus the maps $\alpha_a$
and $\alpha_{a^{-1}}$ are automorphisms of $F$ which
are inverse to each other.

Finally, for each $a,b \in L$ that are adjacent in \gm\ and each
$d \in D$ and $d_t \in T_d$, the equivalence in $(*)$ shows that
the condition $d_{tab} \in T_d$ is equivalent to the conjunction of
the conditions $d_{ta} \in T_d$ and $d_{tb} \in T_d$, so we have
\[ \alpha_b(\alpha_a(d_t)) =
% \left( \begin{cases} d_{ta}, & ta \in T_d \\
%                      d_t,  & ta \not \in T_d
%        \end{cases} \right) ^b =
 \begin{cases} d_{tab}, & d_{ta} \in T_d, \ d_{tab} \in T_d \\
               d_{ta},  & d_{ta} \in T_d, \ d_{tab} \not \in T_d\\
               d_{tb},  & d_{ta} \not \in T_d, \ d_{tb} \in T_d\\
               d_t,     & d_{ta} \not \in T_d, \ d_{tb} \not \in T_d
 \end{cases} =
 \begin{cases} d_{tab}, & d_{ta} \in T_d, \ d_{tb} \in T_d \\
               d_{ta},  & d_{ta} \in T_d, \ d_{tb} \not \in T_d\\
               d_{tb},  & d_{ta} \not \in T_d, \ d_{tb} \in T_d\\
               d_t,     & d_{ta} \not \in T_d, \ d_{tb} \not \in T_d
 \end{cases}.
\]
Therefore, by symmetry, $\alpha_b(\alpha_a(d_t)) =
\alpha_a(\alpha_b(d_t))$. Thus $\alpha_a\alpha_b =
\alpha_b\alpha_a$ whenever $a$ and $b$ are adjacent in $\Gamma$,
which implies that Equation~(\ref{AGL action}) defines a
homomorphism $\alpha:A\g_L \ra Aut(F)$, given by $a \mapsto
\alpha_a$, and an action of $A\Gamma_L$ on the free group $F$.

Let $G := F \rtimes A\Gamma_L$ be the semidirect product defined
by this action. Next we show that $G \cong A\Gamma$. A
presentation for $G$ is given by
\[ G = \langle L \cup (\cup_{d \in D} T_d) \ | \
    R_L \cup (\cup_{d \in D} R_d) \rangle, \]
where
$R_L$ is the set of commutation
relations defining $A\Gamma_L$ (induced by
the edges of $\Gamma_L$) and for each $d \in D$,
\[ R_d := \{ \ d_t^a = d_{ta} \ | \ d_{t} \in T_d, \ d_{ta} \in T_d \ \}
    \cup \{ \ d_t^a = d_t \ | \ d_{t} \in T_d, \ d_{ta} \not\in T_d \ \}.
\]
Next apply Tietze transformations to simplify this presentation.
Given an element $d_t \in T_d$, let $\eta_t$ be a geodesic representative
of $t$.  For each prefix $u$ of
$\eta_t$, then $d_{\ov{u}} \in T_d$ as well, so the
relations of the type
$d_t^a = d_{ta}$ in $R_d$ can be used to show that
$d_t = d_1^{\eta_t}$ in $G$.  For any other geodesic representative
$w$ of $t \in A\g_L$,
the relation $d_1^{\eta_t}=d_1^w$ is a consequence of the
relations in $R_L$.
If we denote $d=d_1$ for $d \in D$, then the presentation of
$G$ is Tietze
equivalent to
\[  \langle L \cup D \ | \
    R_L \cup (\cup_{d \in D} R'_d) \rangle, \]
where
\[ R'_d := \{ \ d^{ta} = d^t \ | \ d_{t} \in T_d, \ d_{ta} \not\in T_d \ \}.
\]
Note that the relation $d^a=d$ occurs in $R'_d$ if $d_a \not\in T_d$,
and $d_a \not\in T_d$ if and only if $a$ is adjacent to $d$
in $\Gamma$.
Thus the relations in $R'_d$ include all the
defining relations in $A\Gamma$ involving $d$.
For each relation $d^{ta} = d^t$ in $R'_d$ with $t$ a nontrivial
element of $A\g_L$, we have $d_{ta} \not\in T_d$, which implies by
$(*)$ that $a$ is adjacent to $d$ and to
all of the symbols in any geodesic for $t$.
This shows that the relation $d^{ta}=d^t$ is a consequence of the
relation $d^a=d$ and the relations in $R_L$. Thus the presentation
for $G$ is
Tietze equivalent to
\[  \langle L \cup D \ | \
    R_L \cup (\cup_{d \in D} R''_d) \rangle, \]
where
\[ R''_d := \{ \ d^a = d \ | \ d_a \not\in T_d \ \}
  = \{ \ d^a = d \ | \ a \text{ is adjacent to }d
           \text{ in } \Gamma \ \}, \]
which is exactly the defining presentation of $A\Gamma$.

Therefore $A\g \cong G = F \rtimes A\Gamma_L$.  By induction
$A\g_L$ has a \pf\ tower of length $\chr(\Gamma)-1$, so the proof
of Proposition~\ref{propn:quotient} completes the proof that $A\g$
has a \pf\ tower of length $\chr(\g)$ and hence $\pfl(A\g) \le
\chr(\g)$.

Next consider the lower bound on the \pf\ length. Let $m=\clq(\g)$
and let $\widetilde \g$ be a clique of \gm\  with $m$ vertices.
Then $\widetilde \g$ is a complete graph, and the subgroup
$A\widetilde \g$ corresponding to $\widetilde \g$ is isomorphic to
$\Z^m$. Lemma \ref{lem:z3} says that $m = \pfl(A\widetilde \g)$,
and Proposition~\ref{propn:subgroup} shows that $\pfl(A\widetilde
\g) \le \pfl(A\g)$.
\end{proof}

%%%%%%%%%%%%%%%%%%%%%%%%%%%%%%%%%%%%%%%%%%%%%%%%%%%%%%%%%%%%%%%%%%%%
\section{Poly-freeness of length 2}\label{sec:niceimp2}
%%%%%%%%%%%%%%%%%%%%%%%%%%%%%%%%%%%%%%%%%%%%%%%%%%%%%%%%%%%%%%%%%%%%

In this section we prove converse of the main result in
Section~\ref{sec:pfimp}, that every graph with the \nice\  induces
a \pf\ \rg\ of length 2.
Together with the main result in Section~\ref{sec:pfimp}, this
completes the proof of Theorem \ref{thm:nicergispf2}.

Given a graph \gm\ together with a
corresponding set $D$ for which \gm\ has the \nice, let
$L:=V(\g)-D$  and let $\g_L$ be the full subgraph of \gm\ induced
by $L$.  Denote by $F(C)$ the free group on the set of connected
components of $\g_L$. If we define a homomorphism $\phi:A\g \ra
F(C)$ by $\phi(d)=1$ for $d \in D$ and $\phi(y)=c$ whenever $y \in
L$ and $c$ is the generator of $F(C)$ corresponding to the
component of $\g_L$ containing $y$, then $D$ is the set of dead
vertices and $\g_L$ is the living subgraph associated to this
homomorphism.  In
the proof below, our approach follows the same lines as the proof
of Theorem~\ref{thm:rgispf}. We construct a free group $F$
(isomorphic to $\ker(\phi)$) and an action of $F(C)$ on $F$, in
order to show explicitly that $A\g$ is isomorphic to a semidirect
product $F \rtimes F(C)$ of two free groups.

\begin{propn}\label{propn:nicergimppf2}
If \gm\  is a graph with the \nice, then
the \rg\  $A\g$ is \pf\  of length 2.
\end{propn}

\begin{proof}
Fix a set $D$ of independent vertices in $\Gamma$ such that every
cycle in $\Gamma$ meets $D$ at least twice.  Let $L$ be the
complementary set of vertices in $\Gamma$ and let $\g_L$
be the full subgraph of \gm\ induced by $L$.
Each of the connected components of the graph
$\g_L$ is a tree.
Select one vertex from each component of $\g_L$,
and denote the set of these vertices by $C$.
Define
$F(C):=A\Gamma_C$ to be the subgroup of $A\Gamma$ corresponding
to the subgraph $\g_C$ induced by $C$.
Since $\g_C$ is totally disconnected,
the group $F(C)$ is also the free group on $C$.

For every vertex $y \in L$, there exists a unique element $c \in
C$ such that $y$ and $c$ are in the same component of $\Gamma_L$,
and since this component is a tree, there exists a unique vertex
path $(y^{(n)},  \dots, y^{(2)}, y^{(1)}, c)$ connecting
$y=y^{(n)}$ and $c$ that lies inside the component of $c$ and is
of minimal length. We call $c$ the {\it component representative
of} $y$ and denote it by $r_y$. For each $d \in D$, let $N_d$
denote the set of vertices in \gm\ adjacent to (i.e. neighbors of)
$d$ and let $\rn_d$ be the set of component representatives of the
vertices in $N_d$. For every element $c$ in $\rn_d$ there exists a
vertex $y \in N_d$ that is contained in the connected component
containing $c$; if $y \neq c$, denote this neighbor of $d$ by
$x(d,c)$. For example, for the graph in Figure~\ref{doubly},
$r_{y_{i,j}} = c_j$, for all $i$ and $j$. The element $x(d,c)$ is
defined only in the following two cases: $x(d_1,c_2)=y_{4,2}$ and
$x(d_2,c_2)=y_{3,2}$.

In the following, the normal form of an
element $t \in F(C)$ refers to the freely reduced word over $C^{\pm 1}$
corresponding to $t$.
Define $X:=L-C$.
For each $x \in X$, let $\x$ be a copy of $x$ and
let $\X = \{ \x | x \in X \}$ be the set of such copies.
For each $\x \in \X$, define a set of symbols
\begin{eqnarray*}
T_\x \ \ :=\ \  \{\ \x_t & | & t\in F(C), \text{ the normal form of }
        t \text{ does not} \\
  & & \text{start with a letter in } \{r_x^{\pm 1}\} \ \}.
\end{eqnarray*}
For each $d \in D$, define a set of symbols
\begin{eqnarray*}
T_d \ \ :=\ \  \{\ d_t  &|&  t\in F(C), \text{ the normal form of }
        t \text{ does not}\\
  & & \text{start with a letter in } \rn_d^{\pm 1} \ \}.
\end{eqnarray*}
For any $z \in \X \cup D$, let $F(T_z)$ be the free group over $T_z$ and
let $F$ be the free group
\[ F := (*_{\x \in \X} F(T_\x)) * (*_{d \in D} F(T_d)) . \]

Given any $c \in C$, define an endomorphism $\alpha_c$ of
the free group $F$ by defining $\alpha_c$ on the generators of $F$
as
\begin{equation}\label{FAx action}
\alpha_c(\x_t) := \begin{cases} \x_{tc}, &
              t \neq 1 \text{ or } c \neq r_x \\
  [(\x^{(n-1)}_1)^{-1} \x^{(n)}_1]
    \cdots   [(\x^{(2)}_1)^{-1} \x^{(3)}_1]
    [(\x^{(1)}_1)^{-1} \x^{(2)}_1] \x^{(1)}_1, & t=1 \text{ and }  c=r_x
       \end{cases}
\end{equation}
for $\x \in \X$ and $t \in T_\x$,  where $(x^{(n)}, \dots,
x^{(2)}, x^{(1)}, c)$ is the minimal length path from $x$ to $c$
inside the component of $c$, and
\begin{equation}\label{FAd action}
\alpha_c(d_t) := \begin{cases} d_{tc}, &
          t\neq 1 \text{ or } c \not\in \rn_d \\
        d_1, & t=1  \text{ and }  c \in \rn_d \cap N_d \\
     d_1^{\x(d,c)_1},  &   t=1  \text{ and } c \in \rn_d - N_d
       \end{cases}
\end{equation}
for $d \in D$ and $t \in T_d$.

Similarly, for $c \in C$, we also define an endomorphism
$\alpha_{c^{-1}}:F \ra F$ by
\begin{equation}\label{FAx- action}
\alpha_{c^{-1}}(\x_t) =
  \begin{cases}
    \x_{tc^{-1}}, & t \neq 1 \text{ or } c \neq r_x \\
    \x^{(1)}_1[\x^{(2)}_1(\x^{(1)}_1)^{-1}][\x^{(3)}_1(\x^{(2)}_1)^{-1}] \cdots
         [\x^{(n)}_1 (\x^{(n-1)}_1)^{-1}], &
    t=1  \text{ and }  c=r_x
  \end{cases},
\end{equation}
for $\x \in \X$ and $t \in T_\x$, where $(x^{(n)},  \dots,
x^{(2)}, x^{(1)}, c)$ is the minimal length path from $x$ to $c$
inside the component of $c$, and
\begin{equation}\label{FAd- action}
\alpha_{c^{-1}}(d_t) =
  \begin{cases}
    d_{tc^{-1}}, & t\neq 1 \text{ or } c \not\in \rn_d \\
    d_1, & t=1 \text{ and } c \in \rn_d \cap N_d, \\
    d_1^{[\alpha_{c^{-1}}(\x(d,c)_1)]^{-1}},  & t=1 \text{ and }
            c \in \rn_d - N_d
  \end{cases},
\end{equation}
for $d \in D$ and  $t \in T_d$.

As in the proof of Theorem~\ref{thm:rgispf}, it is straightforward to
check that the composition of the endomorphisms $\alpha_c$ and
$\alpha_{c^{-1}}$ in either order is equal to the identity on the
generating set of $F$. (The claim easily follows from
$\alpha_c(x^{(i)}_1(x^{(i-1)}_1)^{-1}) =
(x^{(i-1)}_1)^{-1}x^{(i)}_1$ and
$\alpha_{c^{-1}}((x^{(i-1)}_1)^{-1}x^{(i)}_1) =
x^{(i)}_1(x^{(i-1)}_1)^{-1}$, which hold for any adjacent pair of vertices
$x^{(i)}$ and $x^{(i-1)}$ in the component of $c$ at distance $i$
and $i-1$ from $c$, respectively, as measured within the component
of $c$.) Therefore $\alpha_c$ and $\alpha_{c^{-1}}$ are mutually
inverse automorphisms of $F$.

Since the group $F(C)$ is free, the map $c \mapsto \alpha_c$ can
be extended to a homomorphism $\alpha:F(C) \to \Aut(F)$. Therefore
(\ref{FAx action}) and (\ref{FAd action}) define an action of
$F(C)$ on $F$.

%We show now that $A\Gamma \cong F \rtimes F(A)$, which means that the
%poly-free length $A\Gamma$ is 2.
Let $G :=F \rtimes F(C)$ be the associated semidirect product. A
presentation for $G$ can be given by
\[ G = \langle C \cup (\cup_{\x \in \X} T_\x) \cup (\cup_{d \in D} T_d)
   \ | \
    (\cup_{\x \in \X} R_\x) \cup (\cup_{d \in D} R_d) \rangle, \]
where, for $d \in D$ and $\x \in \X$,
\begin{align*}
 R_\x := & \{ \ \x_t^c = \x_{tc} \ |
  \ c \in C, \ \x_t \in T_\x, \ t \neq 1 \text{ or } c \neq r_x \ \}
        \ \cup \\
      & \{ \ \x_1^c =   [(\x^{(n-1)}_1)^{-1} \x^{(n)}_1]
    \cdots   [(\x^{(2)}_1)^{-1} \x^{(3)}_1]
    [(\x^{(1)}_1)^{-1} \x^{(2)}_1] \x^{(1)}_1 \ | \   c=r_x \ \},
\end{align*}
where $(x^{(n)}, \dots, x^{(2)}, x^{(1)}, c)$ is the minimal
length path from $x$ to $c$ inside the component of $c$ and
\begin{align*}
 R_d := & \{ \ d_t^c = d_{tc} \ |
  \ c \in C,\ d_t \in T_d, \ t\neq 1 \text{ or } c \not\in \rn_d \ \}
  \ \cup \\
       & \{ \ d_1^c = d_1 \ | \ c \in \rn_d \cap N_d \ \} \ \cup \\
       & \{ \ d_1^c = d_1^{\x(d,c)_1} \ | c \in \rn_d - N_d \ \}\ .
\end{align*}
Next apply Tietze transformations to simplify this
presentation.

First note that if $\x_t \in T_\x$, then every prefix $u$ of the
normal form of $t$ does not start with a letter in $\{r_x^{\pm
1}\}$, so $\x_{\ov{u}} \in T_\x$ as well, and similarly for $d_t
\in T_d$. Using this fact and repeatedly applying the relations of
the type $\x_t^c = \x_{tc}$ and $d_t^c= d_{tc}$ in $R_\x$ and
$R_d$, respectively, shows that in $G$ we have $\x_t = \x_1^t$ for
all $t \in T_\x$ and $d_t = d_1^t$ for all $t \in T_d$. If we
denote $\x=\x_1$ and $d=d_1$ for each $\x \in \X$ and $d \in D$,
then the presentation for $G$ is Tietze equivalent to
\[  \langle C \cup \X \cup D \ | \
    (\cup_{\x \in \X} R'_\x) \cup (\cup_{d \in D} R'_d) \rangle, \]
where
\[ R'_\x := \{ \ \x^{r_x} =
[(\x^{(n-1)})^{-1} \x^{(n)}]
    \cdots   [(\x^{(2)})^{-1} \x^{(3)}]
    [(\x^{(1)})^{-1} \x^{(2)}] \x^{(1)} \ \} \]
and
\[ R'_d := \{ \ d^c = d \ | \ c \in \rn_d \cap N_d \ \} \ \cup
           \{ \ d^c = d^{\x(d,c)} \ | c \in \rn_d - N_d \ \}. \]

Second, for every $\x \in \X $ introduce a single new generator $x$ and a
relation $\x = x^{-1}r_x$ in the presentation for
$G$. Use these new relations to eliminate
the generators $\x \in \X$ from the above presentation and obtain a
Tietze equivalent presentation
\[  \langle C \cup X \cup D \ | \
    (\cup_{x \in X} R''_x) \cup (\cup_{d \in D} R''_d) \rangle, \]
where
\[ R''_x := \{ \ x^{-1}r_x = [x^{(n-1)} (x^{(n)})^{-1}]
             \cdots [x^{(2)} (x^{(3)})^{-1}] [x^{(1)}(x^{(2)})^{-1}]r_x(x^{(1)})^{-1} \ \}\]
 and
\[ R''_d := \{ \ d^c = d \ | \ c \in \rn_d \cap N_d \ \} \ \cup
          \{ \ d = d^{x(d,c)} \ | c \in \rn_d - N_d \ \}. \]

The relations in $R''_d$ say that each $d \in D$ commutes with all
$c \in C$ and $x \in X$ that are its neighbors in $\Gamma$, just
as in the standard presentation of $A\Gamma$. For each $c \in C$
and each vertex $x$ adjacent to $c$ in $\g_L$, we have $c=r_x$ and
$x=x^{(1)}$.  The corresponding relation in $R''_x$ is $x^{-1}r_x
= r_xx^{-1}$, which implies that $x$ and $c$ commute. When the
path length in $\g_L$ from $x$ to $c$ is 2, with a vertex path
$(x^{(2)},x^{(1)},c)$ from $x=x^{(2)}$ to $c=r_x$, the
corresponding relation is $x^{-1}r_x = x^{(1)} x^{-1} r_x
(x^{(1)})^{-1}$. Since $x^{(1)}$ and $r_x$ commute, this implies
that $x$ and $x^{(1)}$ commute. Continuing in the same fashion we
see that all such relations together imply that each generator $x$
whose distance to $r_x$ in $\g_L$ is $n$ commutes with the
generator that is the neighbor of $x$ on the unique length minimal
path from $x$ to $c=r_x$ inside the component of $r_x$. Thus, the
standard relations in $A\Gamma$ can be recreated from the
relations in $R''_x$ and $R''_d$. Conversely, each relation
$x^{-1}r_x = [x^{(n-1)} (x^{(n)})^{-1}]
             \cdots [x^{(2)} (x^{(3)})^{-1}]
                [x^{(1)}(x^{(2)})^{-1}]r_x(x^{(1)})^{-1}$
in $R''_x$ is a corollary of the defining relations in $A\Gamma$.
Thus the last presentation above is Tietze equivalent to the
standard presentation of $A\Gamma$.

Therefore $A\Gamma \cong  G = F \rtimes F(C)$ has \pf\  length at
most 2. The \nice\ implies that \gm\ is not totally disconnected,
so the group $A\g$ contains a $\Z^2$ subgroup and cannot be free.
Thus the \pf\ length of $A\g$ is exactly 2.
\end{proof}

The free group automorphisms $\alpha_a$ constructed in the proof
of Theorem~\ref{thm:rgispf}  permute the basis elements of the
free group.  Although the free group automorphism $\alpha_c$ of
$F$ in the proof above does not have the same property, the
automorphism of $F^{ab}=F/[F,F]$ induced by $\alpha_c$ permutes
the basis elements of this free abelian group.

We conclude with a fully worked example illustrating a
length 2 \pf\ structure
of a right-angled Artin group defined by
a graph with the doubly breakable cycle property,
following the proof of Theorem B.
\begin{figure}[!hbt]
\begin{center}
\includegraphics{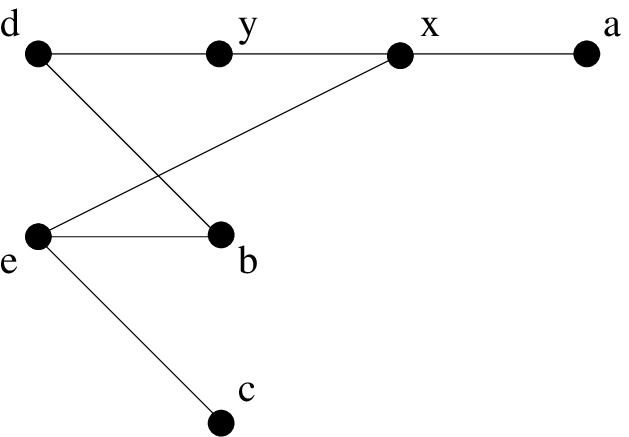}
\end{center}
\caption{Graph with the \nice} \label{example}
\end{figure}
Let $\Gamma$ be the graph in Figure~\ref{example}.
Set $D:=\{d,e\}$. The living subgraph has 3 components and
the chosen representatives are the elements of $C=\{a,b,c\}$.
The only vertices in $X=L-C$ are $x$ and $y$. Denote
$F_q:=F(a,b,c)$. We have
\begin{align*}
 T_d &= \{d_1\} \cup \{ \ d_t  \mid t \in F_q,
  \text{ the normal form of } t \text{ starts with } c^{\pm1} \ \}, \\
 T_e &= \{e_1\},
\end{align*}
and, for $\z \in\{\x,\y\}$,
\[ T_\z = \{ \z_1 \} \cup \{ \ \z_t  \mid
  t \in F_q, \text{ the normal form of } t \text{ starts with }
  b^{\pm1} \text{ or } c^{\pm1} \ \}.
\]
Denote $F_k:=F(T_d \cup T_e \cup T_\x \cup T_\y)$.
Then $A\Gamma=F_k \rtimes F_q$, where the action of
$F_q$ on $F_k$ is given by the following table.
\begin{center}
\begin{tabular}{C|CCCCC}
 & a && b && c  \\[3mm]
\hline
&&&&&\\
d_1 & d_1^{\y_1} && d_1 && d_c \\[3mm]
d_t & d_{ta} && d_{tb} && d_{tc} \\[3mm]
e_1 & e_1^{\x_1} && e_1 && e_1 \\[3mm]
\x_1 & \x_1 && \x_b && \x_c \\[3mm]
\x_t & \x_{ta} && \x_{tb} && \x_{tc} \\[3mm]
\y_1 & (\x_1)^{-1}\y_1\x_1 && \y_b && \y_c\\[3mm]
\y_t & \y_{ta} && \y_{tb} && \y_{tc}
\end{tabular}
\end{center}
In this table, the entry in the row labeled on the
left by the letter $\sigma$ and column labeled above by $\tau$ is
the conjugate $\sigma^\tau$.
In the leftmost column of the table, the $d_t$, $\x_t$,
and $\y_t$ entries range over all symbols in
$T_d \setminus \{d_1\}$,
$T_{\x} \setminus \{\x_1\}$, and
$T_{\y} \setminus \{\y_1\}$,
respectively.

\bibliography{refs_pfrg}
\bibliographystyle{plain}

\end{document}